\documentclass[letterpaper, 10 pt, conference]{ieeeconf}  

\IEEEoverridecommandlockouts                              
\usepackage{graphics} 
\usepackage{epsfig} 
\usepackage{times} 
\usepackage{amsmath} 
\usepackage{amssymb}  
\usepackage{comment}
\usepackage{xcolor}
\usepackage{float}
\usepackage{comment}
\usepackage{enumerate}
\usepackage{booktabs, siunitx, makecell}
\usepackage{algorithm}
\usepackage{algpseudocode}
\algrenewcommand\algorithmicindent{1.0em}%

\usepackage{mathtools}

\newtheorem{theorem}{Theorem}
\newtheorem{lemma}{Lemma}%
\newtheorem{remark}{Remark}%
\newtheorem{definition}{Definition}
\newenvironment{pf}
{\noindent\textit{Proof.}\hspace{1.2mm}}
{\hspace{\fill}$\square$}

\DeclareMathOperator*{\diag}{diag}

\newcommand*{\minimize}[1]{\operatorname*{min}_{#1}\quad} 
\newcommand*{\maximize}[1]{\operatorname*{max}_{#1}\quad} 
\newcommand{\subjectto}{\textup{s.t.}\quad}

\newcommand{\redx}[1]{{\color{red}#1}}

\title{\LARGE \bf
Imitation Learning with Safety and $L^{2}$ Stability Certificates for Boundary Control of Reaction–Diffusion PDEs}

\author{Paulo Henrique F. Biazetto$^{1}$, Mirko Fiacchini$^{2}$, Christophe Prieur$^{2}$, and Gustavo Artur de Andrade$^{1}$
\thanks{This work was partially supported by CAPES under grants 88887.629803/2021-00  and 88881.878833/2023-01 (SticAmSud).}%
\thanks{$^{1}$Paulo Henrique F. Biazetto and Gustavo A. de Andrade
are with the Federal University of Santa Catarina, Department of Automation and Systems, Rua Delfino Conti, 88040-370, Florianópolis, Brazil. {\tt\small paulo.biazetto@posgrad.ufsc.br, gustavo.artur@ufsc.br}}
\thanks{$^{2}$Mirko Fiacchini and Christophe Prieur are with the Université Grenoble Alpes, CNRS, Grenoble-INP, GIPSA-lab, F-38000, Grenoble, France. {\tt\small mirko.fiacchini@gipsa-lab.fr, christophe.prieur@gipsa-lab.fr}}%
}

\begin{document}

\maketitle
\thispagestyle{empty}
\pagestyle{empty}

\begin{abstract}
This paper proposes an imitation learning (IL) framework for synthesizing neural network (NN) controllers that achieve boundary stabilization of systems governed by reaction–diffusion partial differential equations (PDEs). The plant is assumed to be actuated through a Dirichlet boundary condition and subject to a Neumann condition on the unactuated side. The design is based on a finite-dimensional truncated model that captures the unstable dynamics of the original infinite-dimensional system, which is obtained via spectral decomposition. Convex stability and safety conditions are then derived for this truncated model by combining Lyapunov theory with local quadratic constraints (QC), which bound the nonlinear activation functions of the NN and guarantee robustness to model truncation, thus addressing the spillover problem. These conditions are integrated into the IL process to jointly minimize the imitation loss and maximize the volume of the certified region of attraction (ROA). The proposed framework is validated on an unstable reaction–diffusion PDE, demonstrating that the resulting NN controller efficiently reproduces the expert policy while ensuring formal stability guarantees.
\end{abstract}


\section{Introduction}

Many engineering systems are inherently governed by partial differential equations (PDEs), whose distributed-parameter nature poses major challenges for achieving performance, safety, and efficiency \cite{morris2020control}. When accurate models are available, methods such as backstepping \cite{krstic2008}, semigroup and operator-theoretic techniques \cite{curtain2020}, Lyapunov and energy approaches \cite{bastin2016}, and optimal control formulations \cite{li1995optimal} have shown remarkable success in stabilizing infinite-dimensional systems with rigorous guarantees. However, obtaining and calibrating such models often demands costly and time-consuming identification and validation procedures that remain sensitive to modeling uncertainties.

Data-driven methods based on neural networks (NNs) have recently emerged as powerful alternatives. These approaches approximate the dynamics of PDE systems by identifying governing parameters \cite{tanyu2023deep} or by learning the underlying solution operators \cite{luo2025physics}. A more direct paradigm learns the control policy itself, for example, through reinforcement learning \cite{mattioni2023enhancing} or operator-learning techniques \cite{bhan2023neural,zhang2024operator}. While effective in practice, most analyses still provide only a posteriori verification of robustness and stability after training.

A promising approach to bridge this gap is the integration of robust control theory within the learning architecture. In finite-dimensional settings, this has enabled tractable formulations for controller synthesis endowed with formal Lyapunov-based stability certificates. Recent results demonstrate that imitation learning (IL) can be cast as a constrained optimization problem that simultaneously minimizes the imitation loss and maximizes a certified region of attraction (ROA) \cite{yin2021imitation}, while multi-objective extensions introduce tunable trade-offs between tracking fidelity and prescribed stability margins \cite{pinguet2023data}. These frameworks rely on quadratic constraint (QC) formulations to bound neural network (NN) nonlinearities and ensure closed-loop stability during training \cite{zare2024survey}. However, extending such guarantees to systems governed by PDEs remains a major theoretical challenge. The work in \cite{de2023imitation} provided an initial IL formulation for NN-based controllers of linear PDEs, yet key issues remain open—particularly, the robustness of the Lyapunov certificate to finite-dimensional projections and the explicit enforcement of exponential stability rates within the synthesis procedure.


Motivated by these challenges, this work proposes an IL framework for the boundary stabilization of reaction–diffusion PDEs. The plant is actuated through a Dirichlet boundary condition and subject to a Neumann condition on the unactuated side. The methodology begins with a spectral decomposition of the PDE, grounded in Sturm–Liouville theory, to obtain a finite-dimensional model that captures the dominant dynamics of the original system. An expert model predictive controller (MPC) is then designed for this truncated model to generate optimal control trajectories that respect physical and operational constraints. Extending prior work, our framework synthesizes a continuous-time NN controller, learned via the IL framework, that emulates the expert policy while incorporating formal stability guarantees directly from this discrete-time training data. Convex stability and safety conditions are formulated as linear matrix inequalities (LMIs) by combining Lyapunov theory with QCs. Crucially, these LMIs not only bound the NN nonlinearities and ensure robustness to model truncation, but also enforce a desired exponential decay rate (a stability margin). This rigorous treatment of the model reduction mitigates the spillover effect and guarantees stabilization of the infinite-dimensional system in the $L^{2}$-norm. The derived LMIs are embedded in the IL process to jointly minimize the imitation loss and maximize the certified ROA. Simulation results validate the effectiveness of the proposed framework.

The remainder of this paper is organized as follows. Section 2 reviews the mathematical preliminaries on Sturm-Liouville theory. Section 3 formally states the control problem and describes the finite-dimensional modeling approach. Section 4 details our main contribution: the IL framework for synthesizing NN controllers with certified stability. Section 5 presents numerical results to validate the proposed methodology. Finally, Section 6 provides the conclusion.

\section{Overview of Sturm-Liouville Theory}
We begin by recalling the mathematical framework for our analysis. Let $L^{2}(0,1)$ be the Hilbert space of square-integrable functions on the interval $[0,1]$, equipped with the inner product $\langle f,g \rangle_{L^{2}(0,1)} = \int_{0}^{1} f(x)g(x) dx$.

\medskip

\begin{definition}[Sturm-Liouville Operator]\label{def:sturm-liouville}
A Sturm-Liouville operator $\mathcal{A}: D(\mathcal{A}) \subset L^{2}(0,1) \rightarrow L^{2}(0,1)$ is a second-order differential operator of the form
\begin{equation}
    \mathcal{A}f = -(p(x)f')', \label{eq:SLoperator}
\end{equation}
\noindent for all $f$ in the domain
\begin{align*}
    D(\mathcal{A}) = \{f \in L^{2}(0,1) :\,\, & \cos(\theta_{0})f(0) - \sin(\theta_{0})f'(0) = 0, \\
    & \cos(\theta_{1})f(1) + \sin(\theta_{1})f'(1) = 0\},
\end{align*}
\noindent where $p \in {C}^1([0,1])$ with $p(x)>0$, and the domain $D(\mathcal{A})$ is a subspace of $L^2(0,1)$ defined by a set of homogeneous boundary conditions, where $\theta_{0},\theta_{1} \in [0,\pi)$.
\end{definition} 

\medskip


The self-adjointness of $\mathcal{A}$ guarantees the existence of a countable set of real eigenvalues $\{\lambda_{n}\}_{n=1}^{\infty}$ with no finite accumulation point, and a corresponding set of eigenfunctions $\{\Phi_{n}\}_{n=1}^{\infty}$ that form a complete orthonormal basis for $L^2(0,1)$ \cite[Theorem 3.2.8]{curtain2020}. This means any function $f \in L^2(0,1)$ can be uniquely represented by its spectral expansion:
\begin{equation*}
    f(x) = \sum_{n=1}^{\infty} f_n \Phi_{n}(x), \quad \text{where} \quad f_n = \langle f, \Phi_{n} \rangle_{L^{2}(0,1)}.
\end{equation*}

The action of the operator $\mathcal{A}$ on a function $f \in D(\mathcal{A})$ simplifies to scaling each mode by its eigenvalue $\mathcal{A}f = \sum_{n=1}^{\infty} \lambda_{n} f_n \Phi_{n}(x)$.
\section{Problem Statement}
We focus on the stabilization of an unstable one-dimensional linear reaction-diffusion equation. The specific control objective is to steer the system to its zero equilibrium profile, while respecting constraints imposed on the state and the control action. This system serves as a basis example of a PDE where model-based control is well-understood in the literature, providing a clear and solid performance and criterion baseline for our learning NN-based approach with stability guarantees. In the remainder of this section, we will demonstrate how the Sturm-Liouville operator theory is applied to transform the infinite-dimensional boundary control problem into a tractable reduced-order model, providing a tractable foundation for controller synthesis and the IL process.

\subsection{The Reaction-Diffusion System}
We consider the following one-dimensional linear reaction-diffusion system  \cite{lhachemi2020pi}:
\begin{align}
    &z_{t}(t,x) = (p(x)z_{x}(t,x))_{x} + q_{c}z(t,x), \label{eq:pde_main}\\
    &z_{x}(t,0) = 0, \quad z(t,1) = u(t), \label{eq:pde_bc} 
\end{align}
\noindent where $t\in[0,\infty)$ is time, $x\in [0,1]$ is space, $q_{c} \in \mathbb{R}$, and $z(0,x) = z_{0}(x)$ is the initial condition. The control $u(t)$ produces a non-homogeneous condition at the boundary $x = 1$, preventing the direct application of the Sturm-Liouville theory outlined before.

To address this, we introduce a change of variables $w(t,x) = z(t,x) + b(x) u(t)$ to transform the system into an equivalent one with homogeneous boundary conditions, yielding:
\begin{align}
    w_{t}(t,x) &= (p(x)w_{x}(t,x))_{x} + q_{c}w(t,x) \nonumber \\
    &+ a(x)u(t) + b(x) \dot{u}(t), \label{eq:w_pde} \\
    w_x(t,0) &= 0, \quad w(t,1) = 0, \label{eq:w_bcs}
\end{align}
\noindent where $a(x) = - (p(x)b'(x))' - q_{c}b(x)$. This new system for $w(t,x)$ indeed satisfies the homogeneous boundary conditions required for the Sturm-Liouville operator $\mathcal{A}$.



\medskip

\begin{remark} \label{rem:b_conditions}
    The function $b$ must be chosen to satisfy the conditions $b(0)=0$ and $b(1)=-1$ to ensure the boundary conditions in \eqref{eq:w_bcs} are homogeneous.
\end{remark}

\subsection{Spectral Decomposition}
The homogeneous boundary conditions \eqref{eq:w_bcs} allow us to analyze the system's dynamics using a Sturm-Liouville operator. For the clarity of the modal analysis, we define this operator based on the principal differential part of the system, treating the constant reaction term $q_c$ separately.

\medskip

\begin{definition}[System Operator]
The operator $\mathcal{A}$ associated with the system \eqref{eq:w_pde}-\eqref{eq:w_bcs} is given by $\mathcal{A}f = -(pf')'$ with the domain:
\begin{equation}
    D(\mathcal{A}) = \left\{ f \in L^{2}(0,1) : f'(0) = 0 \text{ and } f(1) = 0 \right\}. \label{eq:SLdomain}
\end{equation}
\end{definition}

\medskip

The operator $\mathcal{A}$ is self-adjoint on its domain $D(\mathcal{A})$ and has a set of real, simple eigenvalues and corresponding eigenfunctions which form a complete orthonormal basis for $L^{2}(0,1)$. We can now project the dynamics of the state $w(t,x)$ onto this eigenbasis:
\begin{equation*}  
    w(t,x) = \sum_{n=1}^\infty w_n(t)\Phi_{n}(x).
\end{equation*}

Introducing the projections defined for all $n\geq1$ and all $t\geq0$ in
\begin{align*}
   w_n(t) &= \langle w(t, \cdot), \Phi_{n} \rangle_{L^{2}(0,1)}, \\
   a_{n} &= \langle a, \Phi_{n} \rangle_{L^{2}(0,1)}, \\
   b_{n} &= \langle b, \Phi_{n} \rangle_{L^{2}(0,1)},
\end{align*}
\noindent and the fact that $\mathcal{A}$ is self-adjoint, we get
\begin{align}
    \dot{w}_{n}(t) = (-\lambda_{n} + q_{c})w_{n}(t) + a_{n}u(t) + b_{n}\dot{u}(t). \label{eq:w_decomposed}
\end{align}

Defining the projection coefficients $z_{n}(t) = \langle z(t,\cdot), \Phi_{n} \rangle$ and using the inverse transformation and taking the time derivative, it gives the infinite set of ODEs:
\begin{equation}
    \dot{z}_{n}(t) = (-\lambda_{n}+q_{c}) z_{n}(t) + \beta_{n}u(t), \, \forall n \geq 1 \mbox{ and } t \geq 0, \label{eq:modes}
\end{equation}
\noindent where the coefficients $\beta_{n} = (-\lambda_{n} + q_{c})b_{n} + a_{n}$ determine how the control input influences each mode.

\subsection{Finite-Dimensional Model}
For practical control design, the infinite-dimensional system \eqref{eq:w_pde} is approximated by a finite-dimensional model. We choose an integer $n_{0} \geq 1$ to define a reduced-order model that contains the first $n_{0}$ modes, which are assumed to encompass any unstable or slowly stable dynamics. The feedback controller is designed to stabilize this dominant subsystem.

The finite-dimensional truncated system is:
\begin{align}
    \dot{Z}_{n_{0}}(t) &= A_{n_{0}}Z_{n_{0}}(t) + B_{n_{0}}u(t), \label{eq:sys_trunc}
\end{align}
\noindent where 
\begin{align*}
    Z_{n_{0}}(t) &= \begin{bmatrix} z_{1}(t) & z_{2}(t) & \cdots & z_{n_{0}}(t) \end{bmatrix}^{T}, \\
    A_{n_{0}} &= \diag\left(-\lambda_{1}+q_{c},\, -\lambda_{2}+q_{c},\, \dots ,\, - \lambda_{n_{0}}+q_{c} \right), \\
    B_{n_{0}} &= \begin{bmatrix} \beta_{1} & \beta_{2} & \cdots & \beta_{n_{0}} \end{bmatrix}^{T},
\end{align*}
\noindent and initial condition $Z_{n_{0}}(0) = Z_{0}$.




\section{Imitation Learning (IL)}
This section formally defines the problem of learning a stabilizing NN controller for a class of infinite-dimensional systems. Our primary objective is to move beyond standard data-driven methods by integrating formal stability and robustness guarantees directly into the controller synthesis process. Specifically, we aim to design a continuous-time NN policy that not only mimics the behavior of a high-performance controller but is also endowed with a computable Lyapunov certificate that guarantees the stability of the infinite-dimensional system.

\subsection{Neural Network Controller}
The NN controller $\pi: \mathbb{R}^{n_{0}}\rightarrow \mathbb{R}$ is assumed to be an $\ell$-layer \emph{feedforward neural network} defined as:
\begin{align}
    \omega_{0}(t) &= Z_{n_{0}}(t), \label{eq:nn_controller1}\\
    \omega_{i}(t) &= \phi_{i}(W_{i}\omega_{i-1}(t)+b_{i}),\qquad i\in\{1,\dots,\ell\}, \label{eq:nn_controller2}\\
    u_{\phi}(t)&= W_{\ell + 1} \omega_{\ell}(t) + b_{\ell +1} \label{eq:nn_controller3},
\end{align}
\noindent where $\omega_{0}\in\mathbb{R}^{n_{0}}$ is the input and $u_{\phi} \in \mathbb{R}$ is the output of the NN, and for $i\in\{1,\dots,\ell\}$, $\omega_{i}\in\mathbb{R}^{n_{i}}$ is the output from the $i$-{th} layer. The model of each layer behavior are defined by a weight matrix $W_{i}\in \mathbb{R}^{n_{i}}\times \mathbb{R}^{n_{i-1}}$, a bias vector $b_{i}\in \mathbb{R}^{n_{i}}$, and an activation function $\phi_{i}:\mathbb{R}^{n_{i}}\rightarrow\mathbb{R}^{n_{i}}$, which is applied element-wise, that is,
\begin{align*}
    \phi_{i}(\nu)\triangleq \left( \varphi(\nu_{1}) \quad \dots \quad \varphi(\nu_{n_{i}}) \right),  
\end{align*}
\noindent where $\varphi: \mathbb{R}\rightarrow \mathbb{R}$ is the predefined scalar activation function, such as hyperbolic tangent, sigmoid, rectified linear unit (ReLU). Finally, we assume the activation is identical in all layers, just for the sake of simplifying our analysis.

\subsection{Isolation of NN nonlinearities}
For $i\in\{1,\dots,\ell\}$, let $\nu_{i}\in\mathbb{R}^{n_{i}}$ be the input to the activation function $\phi_{i}$, that is,
\begin{align*}
    \nu_{i}(t) = W_{i}\omega_{i-1}(t) + b_{i},\qquad i\in\{1,\dots,\ell\}.
\end{align*}

Then, the operation of the $i$-{th} layer can be expressed as $\omega_{i}(t) = \phi_{i}(\nu_{i}(t))$. Gathering the inputs and outputs of all activation functions, we get
\begin{align*}
    \nu_{\phi} \triangleq \begin{bmatrix}
    \nu_{1}\\
    \vdots\\
    \nu_{\ell}
    \end{bmatrix} \in \mathbb{R}^{n_{\phi}}, \qquad \mbox{and} \qquad \omega_{\phi} \triangleq \begin{bmatrix}
    \omega_{1}\\
    \vdots\\
    \omega_{\ell}
    \end{bmatrix} \in \mathbb{R}^{n_{\phi}},
\end{align*}
\noindent where $n_{\phi}\triangleq \sum_{i=1}^{\ell}n_{i}$. 

Now, define the vector of stacked activation functions $\phi:\mathbb{R}^{n_{\phi}}\rightarrow \mathbb{R}^{n_{\phi}}$ such that
\begin{align}
    \phi(\nu_{\phi}) \triangleq \begin{bmatrix}
    \phi_{1}(\nu_{1})\\
    \vdots\\
    \phi_{\ell}(\nu_{\ell})
    \end{bmatrix}. \label{eq:vector_nonlinearities}
\end{align}

In this framework, the NN control policy $\pi$ can be rewritten as
\begin{align}
    \begin{bmatrix}
    u_{\phi}(t) \\ \nu_{\phi}(t)
    \end{bmatrix} &= \mathcal{N} \begin{bmatrix}
    Z_{n_{0}}(t) \\ \omega_{\phi}(t) \\ 1
    \end{bmatrix}, \label{eq:nn_framework1} \\
    \omega_{\phi}(t) &= \phi(\nu_{\phi}(t)), \label{eq:nn_framework2}
\end{align}
\noindent where
\begin{align*}
    \mathcal{N} &= \left[\begin{array}{c|cccc|c}
    0 & 0 & 0 & \hdots & W_{\ell+1} & b_{\ell+1}\\\hline
    W_{1} & 0 & \hdots & 0 & 0 & b_{1}\\
    0 & W_{2} & \hdots & 0 & 0 & b_{2}\\
    \vdots & \vdots & \ddots & \vdots & \vdots & \vdots\\
    0 & 0 & \hdots & W_{\ell} & 0 & b_{\ell}
    \end{array}\right], \\
    &= \begin{bmatrix}
    N_{u z} & N_{u\omega} & N_{ub}\\
    N_{\nu z} & N_{\nu \omega} & N_{\nu b}
    \end{bmatrix}.
\end{align*}

The decomposition \eqref{eq:nn_framework1}-\eqref{eq:nn_framework2} isolates the activation function nonlinearities \cite{fazlyab2020safety}, which in turn allows us to analyze the stability of the controlled system by applying QC that bound the activation functions in a similar fashion as in robust control theory.

\subsection{Quadratic Constraints (QC)}
The key idea of applying QC \cite{yin2021stability} is to obtain a new representation of the NN controller \eqref{eq:nn_framework1}-\eqref{eq:nn_framework2} where the nonlinear activation functions are substituted by the constraints they impose on the pre- and post-activation signal. Evidently, any property that can be guaranteed in this framework is also satisfied by the original NN as well. A typical QC is the \emph{local sector condition}. This condition constrains the input-output behavior of each scalar activation function, forcing it to lie within a sector defined by two lines passing through the origin with slopes $m$ and $r$.

In the scalar case, we say that a nonlinear function $\varphi: \mathbb{R}\rightarrow\mathbb{R}$, with $\varphi(0) = 0$, is sector bounded in the sector $[m,\,r]$, with $m \leq r < \infty$, if the following condition holds for all $\nu \in \mathbb{R}:\,(\varphi (\nu ) - m \nu ) ( r \nu - \varphi (\nu ) ) \geq 0$.

This concept can be extended to the vectorial case of nonlinear functions, as the one in \eqref{eq:vector_nonlinearities}. Indeed, for $i\in\{1,\dots, n_{\phi}\}$, assume that the $i$-th activation function $\phi_{i}$ in $\phi$ is sector bounded by $[m_{i},\,r_{i}]$. Then, it is possible to stack these sectors into vectors $m_{\phi},\,r_{\phi}\in\mathbb{R}^{n_{\phi}}$ so that the following condition is satisfied for $\phi$ \cite{tarbouriech2011stability}:
\begin{align}
    (\phi(\nu_{\phi}) - M_{\phi}\nu_{\phi})^{T}(\phi (\nu_{\phi}) - R_{\phi}\nu_{\phi}) \geq 0,
\label{eq:local_sc_matrix}
\end{align}
\noindent with $M_{\phi} = \diag(m_{\phi})$ and $R_{\phi}= \diag(r_{\phi})$.

From \eqref{eq:local_sc_matrix}, it readily follows that 
\begin{align}
    \begin{bmatrix}
    \nu_{\phi}\\
    \omega_{\phi}
    \end{bmatrix}^{T} \begin{bmatrix}
    -2M_{\phi}R_{\phi} & (M_{\phi}+R_{\phi}) \\
    (R_{\phi}+M_{\phi}) & -2I_{n_{\phi}}
    \end{bmatrix} \begin{bmatrix}
    \nu_{\phi} \\
    \omega_{\phi}
    \end{bmatrix} \geq 0. \label{eq:ineq_sector_condition}
\end{align}

For the NN controller \eqref{eq:nn_framework1}-\eqref{eq:nn_framework2} the local sector conditions can be obtained by a numerical procedure given the set of weighting matrices, $W_{i}$, the bias vectors, $b_{i}$, and monotonically non-decreasing activation functions, $\phi_{i}$. Additionally, we will assume an \emph{artificial set of state constraints}, $\mathcal{Z}$, in order to allow us to prove less conservative bounds within the constraints set. This set is defined by
\begin{align*}
    \mathcal{Z} = \left\{ Z_{n_0} \in\mathbb{R}^{n_{0}}: -s\leq S Z_{n_0} \leq s,\, s\geq 0 \right\},
\end{align*}
\noindent where $S\in \mathbb{R}^{n_{z}\times n_{0}}$ and $s\in\mathbb{R}^{n_{z}}$.

Given the bounds on the NN's input layer, the bounds on the inputs and outputs of each subsequent neuron can be computed. This is achieved by propagating the initial state constraints forward through the network, layer by layer, using standard interval arithmetic techniques. This procedure determines the tightest possible range for the pre-activation signals, which we denote by the interval $[\underline{\nu}, \overline{\nu}]$, for all states within the set $\mathcal{Z}$. These bounds are then used to calculate the corresponding sector parameters, as summarized in the following lemma (see \cite{yin2021imitation,de2023imitation} for details).

\medskip
\begin{lemma} \label{lemma:sector_condition}
    Let $\mathcal{Z} = \{Z_{n_{0}} \in\mathbb{R}^{n_{0}}: \underline{Z}\leq Z_{n_{0}} \leq \overline{Z}\}$ be the smallest hypercube that bounds the states. Additionally, consider that the activation functions of the NN defined in \eqref{eq:nn_framework1}-\eqref{eq:nn_framework2} are given by the hyperbolic tangent function. Then, there exist vectors $m_{\phi}$ and $r_{\phi}$ such that the nonlinearity $\phi$ defined in \eqref{eq:vector_nonlinearities} satisfies the local sector constraint \eqref{eq:local_sc_matrix}. 
\end{lemma}

\subsection{Equivalence of the Closed-Loop System} 

We decompose the NN controller into an affine part and a nonlinear activation described by sector bounds. While this enables a convex relaxation, jointly synthesizing the controller and its stability certificate remains non-convex due to the coupling between the network weights $\mathcal{N}$ and the sector limits $[m_\phi, r_\phi]$. To address this, a loop transformation embeds the sector information into a new linear matrix $\tilde{N}$ and normalizes the nonlinearity to a fixed sector $[-1,1]$. This allows us to efficiently redesign the controller and its Lyapunov-based stability certificate using optimization tools. 

The following Lemma formalizes this entire process.

\medskip

\begin{lemma} \label{lemma:loop_transformation}
    Consider the finite-dimensional system \eqref{eq:sys_trunc} in feedback with a NN controller policy $\pi$ defined by the affine transformation \eqref{eq:nn_framework1} and the nonlinearity \eqref{eq:nn_framework2} sector-bounded by $[M_{\phi}, R_{\phi}]$ satisfying the quadratic sector constraint \eqref{eq:local_sc_matrix} in a set of state $\mathcal{Z}$.

    Then, under the loop transformation
    \begin{equation}
        \omega_{\phi}(t) = \tfrac{R_{\phi}-M_{\phi}}{2}\,x_{\phi}(t) 
        + \tfrac{M_{\phi}+R_{\phi}}{2}\,\nu_{\phi}(t),
        \label{eq:loop_transform}
    \end{equation}
    \noindent it is possible to obtain an equivalent controller representation defined by a matrix $\tilde{N}\in\mathbb{R}^{(n_u+n_\phi)\times(n_0+n_\phi)}$ and a normalized nonlinearity $\tilde{\phi}$:
    \begin{align}
        \begin{bmatrix} u_{\phi}(t) \\ \nu_{\phi}(t) \end{bmatrix} &= \tilde{N} \begin{bmatrix} Z_{n_{0}}(t) \\ x_{\phi}(t) \end{bmatrix}, \label{eq:nn_tilde_framework1} \\
        x_{\phi}(t) &= \tilde{\phi}(\nu_{\phi}(t)), \label{eq:nn_tilde_framework2}
    \end{align}

    The nonlinearity $\tilde{\phi}$ is normalized to the sector $[-1_{n_{\phi} \times 1}, 1_{n_{\phi} \times 1}]$ and satisfies, for any diagonal $\Lambda = \diag(\lambda_i)$ with $\lambda_{i} \geq 0,\,\forall i=1,\dots,n_\phi$,
    \begin{equation}
        \begin{bmatrix} \nu_{\phi} \\ x_{\phi} \end{bmatrix}^{\top} \begin{bmatrix} \Lambda & 0 \\ 0 & -\Lambda \end{bmatrix} \begin{bmatrix} \nu_{\phi} \\ x_{\phi} \end{bmatrix} \geq 0, \quad \forall \nu_{\phi} \in [\underline{\nu}, \overline{\nu}]. \label{eq:normalized_sector}
    \end{equation}
\end{lemma}

\medskip

\begin{pf}
Let 
\begin{equation*}
    N = \begin{bmatrix}
    N_{u z} & N_{u\omega}\\
    N_{\nu z} & N_{\nu \omega}
    \end{bmatrix} 
\end{equation*}
\noindent be the block partition of the affine map appearing in \eqref{eq:nn_framework1} which possessing the stabilization in an equlibrium points $(Z_{n_{0}}^{*}, u_{\phi}^{*}, \nu^{*}, \omega_{\phi}^{*})$. This leads to
\begin{align}
    \begin{bmatrix}
    u_{\phi}(t) \\ \nu_{\phi}(t)
    \end{bmatrix} &= N \begin{bmatrix}
    Z_{n_{0}}(t) \\ \omega_{\phi}(t)
    \end{bmatrix}, \label{eq:nn_framework11} \\
    \omega_{\phi}(t) &= \phi(\nu_{\phi}(t)), \label{eq:nn_framework22}
\end{align}


Define the centralization $\alpha \coloneq \frac{M_{\phi} + R_{\phi}}{2}$ and normalization $\beta \coloneq \frac{R_{\phi} - M_{\phi}}{2}$ scalars, the loop transformation \eqref{eq:loop_transform} can be rewritten as
\begin{equation} \label{eq:proof_wphi}
    \omega_{\phi}(t)=\beta\,x_{\phi}(t)+\alpha\,\nu_{\phi}(t).
\end{equation}

By elementary algebra, the normalized nonlinearity is therefore
\begin{equation} \label{eq:proof_tildephi}
    \tilde\phi(\nu_{\phi})=\frac{\phi(\nu_{\phi})-\alpha\nu_{\phi}}{\beta}
    =\frac{2\phi(\nu_{\phi})-(M_{\phi}+R_{\phi})\nu_{\phi}}{R_{\phi}-M_{\phi}},
\end{equation}
\noindent so $\tilde{\phi}(\nu_{\phi}) \in [-1_{n_{\phi} \times 1}, 1_{n_{\phi} \times 1}]$ and \eqref{eq:normalized_sector} holds.

Define the intermediate matrices $C_1 \coloneq N_{u\omega}\beta$, $C_2 \coloneq N_{u\omega}\alpha$, $C_3 \coloneq N_{\nu\omega}\beta$, and $C_4 \coloneq N_{\nu\omega}\alpha$, and substitute \eqref{eq:proof_wphi} into the affine relation \eqref{eq:nn_framework11} to obtain
\begin{align}
    u_{\phi}(t) &= N_{uz} Z_{n_{0}}(t) + C_1 x_{\phi}(t) + C_2 \nu_{\phi}(t), \label{eq:u_n} \\
    \nu_{\phi}(t) &= N_{\nu z} Z_{n_{0}}(t) + C_3 x_{\phi}(t) + C_4 \nu_{\phi}(t). \label{eq:v_phi}
\end{align}

The expression \eqref{eq:v_phi} can be solved for the variable $\nu_{\phi}(t)$:
\begin{align}
    \nu_{\phi}(t) = (I - C_4)^{-1} N_{\nu z} Z_{n_{0}}(t) + (I - C_4)^{-1} C_3 x_{\phi}(t), \label{eq:v_phi2}
\end{align}
\noindent and substituting into \eqref{eq:u_n} yields
\begin{align}
    u_{\phi}(t) &= \left(N_{uz} + C_2 (I - C_4)^{-1} N_{\nu z}\right) Z_{n_{0}}(t) \nonumber\\
    &+ \left(C_1 + C_2 (I - C_4)^{-1} C_3\right) x_{\phi}(t). \label{eq:u_n2}
\end{align}

Therefore, using \eqref{eq:v_phi2} and \eqref{eq:u_n2}, it is possible to obtain the transformed linear affine map
\begin{align}
    \tilde{N} &= \begin{bmatrix}
        N_{uz} + C_{2}(I - C_{4})^{-1}N_{\nu z} & C_{1}+C_{2}(I-C_{4})^{-1}C_{3} \\ (I-C_{4})^{-1}N_{\nu z} & (I-C_{4})^{-1}C_{3}
    \end{bmatrix} \nonumber \\
    &= \begin{bmatrix}
    \tilde{N}_{u z} & \tilde{N}_{u x}\\
    \tilde{N}_{\nu z} & \tilde{N}_{\nu x}
    \end{bmatrix}. \label{eq:nn_normalized}
\end{align}

The transformed linear map $\tilde{N}$ can be described as a function of $N$, defined as $\tilde{N} = \mathcal{F}(N)$. Finally, this normalization on the sector $[-1_{n_{\phi} \times 1}, 1_{n_{\phi} \times 1}]$ transforms the sector condition \eqref{eq:ineq_sector_condition} in \eqref{eq:normalized_sector}. This transformation is well-defined provided that the matrix $(I - C_4)$ is invertible. This yields the stated result and completes the proof.
\end{pf}

\subsection{Lyapunov Stability}
The stability result of system \eqref{eq:modes} under the NN controller \eqref{eq:nn_tilde_framework1}-\eqref{eq:nn_tilde_framework2} is assessed by the following Theorem.

\medskip

\begin{theorem}\label{theo}
Consider the reaction-diffusion system described by \eqref{eq:modes}. Let $p \in C^{1}([0,1])$ with $p>0$, and $q_{c} \in \mathbb{R}$. Let $n_{0} \geq 0$ and $\delta > 0$ be given such that $-\lambda_{n} + q_{c} < -\delta < 0$ for all $n \geq n_{0}+1$. Let $\underline{Z}$, $\overline{Z}$, $a_{\phi}$ and $b_{\phi}$ be given vectors satisfying Lemma \ref{lemma:sector_condition}. Consider a NN controller \eqref{eq:nn_framework1}-\eqref{eq:nn_framework2} with equilibrium point $Z_{n_{0}}^{*} = 0$ and the state constraint set $\mathcal{Z}$, subjected to the transformation in Lemma \ref{lemma:loop_transformation}.

If there exist a positive definite matrix $P\in \mathbb{R}^{n_{0}\times n_{0}}$, vector $\lambda \in \mathbb{R}^{n_{\phi}}$ with $\lambda \geq 0$, $\alpha,\,\gamma > 0$ and $0 < \tau < 1$, such that the following inequalities hold:
\begin{align}
    \tilde{R}_{\nu}^{T} Q_{1}^{*} \tilde{R}_{\nu} + \tilde{R}_{\phi}^{T}Q_{2}\tilde{R}_{\phi} + Q_{3} &\prec 0, \label{eq:lmi} \\
    \Gamma_{n_{0}+1} = \gamma \sum_{n>n_{0}+1} \left( 2\left(-\lambda_{n}+q_{c}+\delta\right)+\frac{1}{\alpha} \right)z_{n}^{2} &< 0, \label{eq:lmi2} \\
    \alpha\gamma\left\|\mathcal{R}_{\beta}\right\|^{2} \tilde{R}_{\nu}^{T}\tilde{R}_{\nu} &\prec 0, \label{eq:lmi3} \\
    \begin{bmatrix} s_{i} & S_{i}^{T} \\ S_{i} & P \end{bmatrix} &\succeq 0, \label{eq:lmi4}
\end{align}
\noindent for $i=1,\dots,n_{0}$ where $\Lambda = \diag(\lambda)$ and
\begin{align*}
   Q_{1}^{*} &= \begin{bmatrix}
   A_{n_{0}}^{T}P + PA_{n_{0}} + 2\delta P & PB_{n_{0}} \\ B_{n_{0}}^{T}P & 0
   \end{bmatrix}, ~     Q_{2} = \begin{bmatrix}
      \Lambda & 0 \\ 0 & -\Lambda
   \end{bmatrix},& \nonumber \\
   Q_{3} &= \begin{bmatrix}
      \tau P & 0 \\ 0 & \tau\Lambda
   \end{bmatrix},     \qquad\qquad\qquad\qquad \tilde{R}_{\nu} = \begin{bmatrix}
      I_{n_{0}} & 0 \\ \tilde{N}_{u z} & \tilde{N}_{u x}
   \end{bmatrix},& \\
   \tilde{R}_{\phi} &= \begin{bmatrix}
      \tilde{N}_{\nu z} & \tilde{N}_{\nu x} \\ 0 & I_{n_{\phi}}
   \end{bmatrix},
\end{align*}
\noindent then the system \eqref{eq:modes} is asymptotically stable at the equilibrium point. Moreover, a the set $\mathcal{E}(P) \coloneqq \{Z_{n_{0}} \in \mathbb{R}^{n_{0}} : Z_{n_{0}}^TPZ_{n_{0}} \leq 1\}$, is an inner-approximation of ROA, i.e., $\mathcal{E}(P) \subseteq \mathcal{R}$ where $\mathcal{R} \coloneqq \{Z_{0} \in \mathcal{Z} \mid \lim_{t\to\infty} \mathcal{G}(t; Z_{0}) = 0\}$ and $\mathcal{G}(t; Z_{0})$ denote the solution to the feedback system at time $t$ from initial condition $Z_{n_{0}}(0) = Z_{0}$.


\end{theorem}

\begin{pf}
Consider the Lyapunov function
\begin{equation*}
    V(t) = V_{1}(t) + V_{2}(t),
\end{equation*}
\noindent where
\begin{align}
    V_{1}(t) &= Z_{n_{0}}^{T}P Z_{n_{0}},\\
    V_{2}(t) &= \gamma \sum_{n\geq n_{0}+1} \langle z,\Phi_{n} \rangle^{2}_{L^{2}(0,1)},
\end{align}
\noindent where $P \in \mathbb{R}^{n_{0}\times n_{0}}$ is a positive definite matrix, and $\gamma > 0$ is a parameter variable.


Computing the time derivative of $V_{1}$ and substituting \eqref{eq:sys_trunc}, we obtain
\begin{align}
    \dot{V}_{1}(t) = \begin{bmatrix}
        Z_{n_{0}}(t) \\ u_{\phi}(t)
    \end{bmatrix}^{T} Q_{1} \begin{bmatrix}
        Z_{n_{0}}(t) \\ u_{\phi}(t)
    \end{bmatrix} < 0, \label{eq:dV1}
\end{align}
\noindent where
\begin{align*}
    Q_{1} = \begin{bmatrix}
    A_{n_{0}}^{T}P + PA_{n_{0}} & PB_{n_{0}} \\ B_{n_{0}}^{T}P & 0 
    \end{bmatrix}.
\end{align*}

Once the sector condition has been normalized and since the system can be shifted to the origin or to some equilibrium point and stabilized, i.e., the state satisfies $Z_{n_0}^* = 0_{n_0 \times 1}$, as a result of the control action of the NN, this leads to the control action $u_{\phi}^* = 0$ to the equilibrium. Hence, based in the results presented in Lemmas \ref{lemma:sector_condition} and \ref{lemma:loop_transformation}, and then substituting \eqref{eq:u_n2} in \eqref{eq:dV1}, it can be rewritten as
\begin{align*}
    \dot{V}_{1}(t) = \begin{bmatrix} Z_{n_{0}}(t) \\ x_{\phi}(t) \end{bmatrix}^{T} \tilde{R}_{\nu}^{T} Q_{1}\tilde{R}_{\nu} \begin{bmatrix}
        Z_{n_{0}}(t) \\ x_{\phi}(t)
    \end{bmatrix} < 0,
\end{align*}

\noindent and substituting \eqref{eq:v_phi2} in the sector condition \eqref{eq:normalized_sector}, it can be rewritten as
\begin{align*}
    \begin{bmatrix} Z_{n_{0}}(t) \\ x_{\phi}(t) \end{bmatrix}^{T} \tilde{R}_{\phi}^{T}Q_{2}\tilde{R}_{\phi} \begin{bmatrix} Z_{n_{0}}(t) \\ x_{\phi}(t) \end{bmatrix} &\geq 0.
\end{align*}

For the second part, $V_{2} = \gamma \sum_{n\geq n_{0}+1} \langle z,\Phi_{n} \rangle^{2}_ {L^{2}(0,1)}$ along the solution \eqref{eq:modes}, yields
\begin{align*}
    \dot{V}_{2}(t) = 2\gamma \sum_{n\geq n_{0}+1} z_{n}\Big( (-\lambda_{n}+q_{c})z_{n} + \beta_{n}u_{\phi} \Big).
\end{align*}

Using Young's inequality gives
\begin{align*}
    2 \sum_{n\geq n_{0}+1} \beta_{n}z_{n}u_{\phi} \leq \frac{1}{\alpha} \sum_{n\geq n_{0}+1} z_{n}^{2} + \alpha\left\|\mathcal{R}_{\beta}\right\|^{2}_{L^{2}(0,1)}u_{\phi}^{2},
\end{align*}
\noindent for any $\alpha > 0$, where the constant $\left\|\mathcal{R}_{\beta}\right\|^{2}_{L^{2}(0,1)} \coloneq \sum_{n \geq n_{0}+1} \beta_n^2$ captures the energy of the high-order modes.

To ensure exponential stability with decay rate $\delta > 0$, we enforce $\dot{V} + 2\delta V < 0$. The derivation, which includes a relaxation parameter $\tau\in(0,1)$ within the LMI conditions, yields:
\begin{multline*}
   \dot{V} + 2\delta V = \begin{bmatrix} Z_{n_{0}}(t) \\ x_{\phi}(t) \end{bmatrix}^{T} \tilde{R}_{\nu}^{T} Q_{1}^{*} \tilde{R}_{\nu} \begin{bmatrix} Z_{n_{0}}(t) \\ x_{\phi}(t) \end{bmatrix} \\
   + \begin{bmatrix} Z_{n_{0}}(t) \\ x_{\phi}(t) \end{bmatrix}^{T} \tilde{R}_{\phi}^{T} Q_{2} \tilde{R}_{\phi} \begin{bmatrix} Z_{n_{0}}(t) \\ x_{\phi}(t) \end{bmatrix} \\
   + \gamma \sum_{n \geq n_{0}+1} \left( 2\left(-\lambda_{n}+q_{c}+\delta\right)+\frac{1}{\alpha} \right)z_{n}^{2} \\
   + \alpha\gamma\left\|\mathcal{R}_{\beta}\right\|^{2} \tilde{R}_{\nu}^{T}\tilde{R}_{\nu} < - Q_{3}.
\end{multline*}

This completes the proof.
\end{pf}



\medskip

\begin{remark}
Feasibility of inequalities~\eqref{eq:lmi}--\eqref{eq:lmi4} for large truncation orders $n_0$ is not guaranteed, since $P$, $\lambda$, and $Q_1^*,Q_2,Q_3$ depend on the finite-dimensional approximations $A_{n_0},B_{n_0},\tilde{R}_\nu,\tilde{R}_\phi$. Uniform bounds, e.g., $\|P(n_0)\|_2 \le C$, $\|\lambda(n_0)\|_2 \le C$ for some $C>0$ independent of $n_0$, are needed to ensure the feasible region does not vanish as $n_0 \to \infty$.
\end{remark}

\medskip

The feasibility problem of Theorem \ref{theo} is not linear due to the presence of bilinearity with the weights of the NN with the Lyapunov matrix and some terms such as $\frac{1}{\alpha}$ involving the decision variables. To handle these nonlinearities, the problem must be split into two optimization problems.

First, to ensure convexity with respect to $P$ and $\Lambda$, the Lyapunov condition \eqref{eq:lmi} is formulated using a predefined mapping $\tilde{N} = \mathcal{F}(N)$, where $\tilde{N}$ is considered known. To integrate this stability requirement into the training procedure (along with the \eqref{eq:lmi4} condition), we now consider $\tilde{N} \in \mathbb{R}^{(n_u + n_\phi) \times (n_0 + n_\phi)}$ as an additional decision variable alongside $P$ and $\Lambda$. The goal is to establish a stability criterion that remains convex. Thus, we can reformulate the problem as,
\begin{align*}
    \tilde{R}_{v}^{T} Q_{1}^{*} \tilde{R}_{v} + \tilde{R}_{\phi}^{T} Q_{2} \tilde{R}_{\phi} + Q_{3} \prec 0, 
\end{align*}
\noindent and applying Schur complements yields the condition
\begin{align*}
    \begin{bmatrix} -\Theta_P & -P B_{n_{0}} \tilde{N}_{ux} & \tilde{N}_{\nu z}^{T} \\ -\tilde{N}_{ux}^{T} B_{n_{0}}^{T} P & (1-\tau)\Lambda & \tilde{N}_{\nu x}^{T} \\ \tilde{N}_{\nu z} & \tilde{N}_{\nu x} & \Lambda^{-1} \end{bmatrix} \succ 0,
\end{align*}
\noindent where $\Theta_P = \big(A_{n_{0}}^{T} P + P A_{n_{0}} + P B_{n_{0}} \tilde{N}_{uz} + \tilde{N}_{uz}^{T} B_{n_{0}}^{T} P \,\allowbreak + \left(2\delta+\tau\right)P \big)$.

Left and right multiplying by $\diag(P^{-1}, \Lambda^{-1}, I_{n_{\phi}})$, and defining the variables $H_{1} = P^{-1} > 0,\, H_{2} = \Lambda^{-1} > 0,\, L_{1} = \tilde{N}_{uz}H_{1},\, L_{2} = \tilde{N}_{ux}H_{2},\, L_{3} = \tilde{N}_{\nu z}H_{1}$ and $L_{4} = \tilde{N}_{\nu x}H_{2}$, we obtain
\begin{align}
    \mbox{LMI}(\bar{H},\,\bar{L}) = \begin{bmatrix} -\Theta_H & -B_{n_{0}}L_{2} & L_{3}^{T} \\ -L_{2}^{T}B_{n_{0}}^{T} & (1-\tau)H_{2} & L_{4}^{T} \\ L_{3} & L_{4} & H_{2} \end{bmatrix} \succ 0,
\end{align}
\noindent where $\Theta_H = \big(H_{1}A_{n_{0}}^{T} + A_{n_{0}}H_{1} + B_{n_{0}}L_{1} + L_{1}^{T}B_{n_{0}}^{T} \,\allowbreak+ \left(2\delta+\tau\right)H_{1}\big)$,  $\bar{H} = \{H_{1},H_{2}\}$ and $\bar{L}=\{L_{1},L_{2},L_{3},L_{4}\}$.

The second optimization problem involves verifying \eqref{eq:lmi2} and \eqref{eq:lmi3} after the learning training process. We can fix the value of $\gamma$ by a line search and use the Schur complement to write \eqref{eq:lmi2} as
\begin{align*}
    \Gamma_{n_{0}+1}^{*} = \begin{bmatrix} 2\gamma(-\lambda_{n_{0}+1} + q_{c} + \delta) & 1 \\ \star & -\alpha \end{bmatrix} \preceq 0.
    \end{align*}

Therefore, as soon as $\gamma$ (or $\alpha$) is fixed, checking the conditions of Theorem \ref{theo} (with the optimal values $(N^{\star},\,\bar{H}^{\star},\,\bar{L}^{\star})$) reduces to check the LMI:
\begin{align}
    \maximize{} & \alpha \label{residual:op}\\
    \subjectto & \dot{V}_{1} + \tilde{R}_{\phi}^{T} Q_{2} \tilde{R}_{\phi} + \alpha\gamma\left\| \mathcal{R}_{\beta} \right\|^{2}\tilde{R}_{\nu}^{T}\tilde{R}_{\nu} \leq -\delta V_{1}, \label{residual:op_2}\\
    & \Gamma_{n_{0}+1}^{*} \preceq 0, \label{residual:op_3}\\
    & \alpha > 0. \label{residual:op_4}
\end{align}

Thus, given a desired exponential decay rate $\delta > 0$, $0 < \tau < 1$, and a number of modes $n \geq n_{0} + 1$ for the system, the sufficient conditions of the previous theorem can be reformulated as an efficiently numerically tractable problem, i.e, the inequalities \eqref{eq:lmi}-\eqref{eq:lmi4} are always feasible for $n_{0}$ large enough.

\subsection{Convexify}

To ensure stability, we train the NN controller by solving a constrained optimization problem that integrates stability conditions directly into the learning process. The problem is formulated as:
\begin{align}
    \minimize{N,H,L} & \quad \eta_1 \mathcal{L}(N) - \eta_2 \log \det(H_1) \label{eq:obj_concise} \\
    \subjectto & \quad \text{LMI}(\bar{H}, \bar{L}) \succ 0, \label{eq:lmi_concise} \\
               & \quad \mathcal{F}(N)H = L, \label{eq:coupling_concise}
\end{align}
\noindent where 
\begin{align*}
    H = \begin{bmatrix}
        H_{1} & 0 \\ 0 & H_{2}
    \end{bmatrix}, \quad L = \begin{bmatrix}
        L_{1} & L_{2} \\ L_{3} & L_{4}
    \end{bmatrix}.
\end{align*}

The objective \eqref{eq:obj_concise} balances the imitation loss $\mathcal{L}(N)$ against the size of the verifiable ROA (maximized via $\log \det(H_1)$), by the $\eta_1, \eta_2 > 0$ weights. The equality constraint \eqref{eq:coupling_concise} couples the NN weights $N$ to the stability certificate $(H, L)$ thought the transformation $\mathcal{F}(N) = LH^{-1} = \tilde{N}$.

This non-convex problem is solved using the Alternating Direction Method of Multipliers (ADMM) \cite{yin2021imitation}, which decouples the learning and stability-enforcing steps. This is achieved by minimizing the augmented Lagrangian:
\begin{multline*}
   \mathcal{L}_a(N, H, L, Y) = \eta_1 \mathcal{L}(N) - \eta_2 \log \det(H_1) \\
   + \operatorname{tr}\left(Y^T (\mathcal{F}(N)H - L)\right) + \frac{\rho}{2} \left\|\mathcal{F}(N)H - L\right\|_F^2,
\end{multline*}
\noindent where $Y$ is the Lagrange multiplier, $\rho > 0$ is a penalty parameter, and $\|\cdot\|_F$ denotes the Frobenius norm. The iterative procedure, which alternates between updating the NN weights $N$ and the certificate variables $(\bar{H}, \bar{L})$, is outlined in Algorithm \ref{alg:safe_il_pde}.

\begin{algorithm}[htpb]
\caption{Safe Imitation Learning for Parabolic PDEs}
\label{alg:safe_il_pde}
\begin{algorithmic}[1]
\Require Expert data, learning rates $\eta_1, \eta_2$, penalty $\rho$, epochs $E$.
\State \textbf{Initialize:} $N^0, \bar{H}^0, \bar{L}^0, Y^0$, iteration $k \gets 0$.
\Repeat
    \For{$e = 1 \to E$} \Comment{Step 1: NN Weights Update}
        \State $N^{k+1} \gets \text{AdamOptimizer}(\nabla_N \mathcal{L}_a(N, H^k, L^k, Y^k))$
    \EndFor
    \State \Comment{Step 2: Stability Certificate Update}
    \State $(H, L)^{k+1} \gets \arg \min_{H, L} \mathcal{L}_a(N^{k+1}, H, L, Y^k)$,
    \Statex $\qquad \qquad \qquad \qquad \text{s.t.} \quad \text{LMI}(\bar{H}, \bar{L}) \succ 0$.
    \State \Comment{Step 3: Lagrange Multiplier Update}
    \State $Y^{k+1} \gets Y^k + \rho (\mathcal{F}(N^{k+1})H^{k+1} - L^{k+1})$
    \State $k \gets k + 1$
    \Statex \Comment{\textit{Step 4: Final Stability Verification}}
    \Statex Given the converged solution $(N^{\star}, \bar{H}^{\star}, \bar{L}^{\star})$, find the largest $\alpha > 0$ s.t. \eqref{residual:op_2}-\eqref{residual:op_4}.
    \State $\alpha^\star \gets \alpha$
\Until{Given $\alpha^\star$, check if \eqref{eq:lmi2} holds.}
\State \Return $(N^\star, \bar{H}^\star, \bar{L}^\star, \alpha^\star)$.
\end{algorithmic}
\end{algorithm}

While the overall problem is non-convex, preventing global optimality guarantees, the ADMM framework ensures that any converged solution provides a provably robust NN controller, since the stability certificate update is formulated as a convex subproblem solved at each iteration.

\section{Results}

To illustrate the results, we consider a PDE system \eqref{eq:pde_main}-\eqref{eq:pde_bc} subject to the homogeneous transformations with the case where the diffusion coefficient is $p(x)=1$. The spectrum of the operator $\mathcal{A}$ can then be determined analytically and it is well-documented (see, e.g., \cite{boyce2021elementary}). The eigenvalues are given by $\lambda_{n} = \left(n - \tfrac{1}{2}\right)^2 \pi^2$, with corresponding eigenfunctions $\Phi_{n}(x) = \sqrt{2} \cos\left(\left(n - \tfrac{1}{2}\right)\pi x\right),\,\forall n\geq1$. A finite-dimensional model suitable for controller design is then derived by truncating this system to its first $n_{0}$ modes, as described by \eqref{eq:sys_trunc}.

Regarding the actuation part of the system, we can choose the function $b(x) = -x^{2}$, according to the required properties derived in Remark \ref{rem:b_conditions}. For the controller, let's define a linear MPC as
\begin{align}
    \minimize{u} & \sum_{k=0}^{M-1} \left( Z_{n_{0},k}^{T}\bar{Q}Z_{n_{0},k} + u_{k}^{T}\bar{R}u_{k} \right) + Z_{n_{0},M}^{T}\bar{P}Z_{n_{0},M} \\
    \subjectto & Z_{n_{0},k+1} = A_{n_{0}}Z_{n_{0},k} + B_{n_{0}}u_{k}, \\
               & Z_{n_{0},0} = Z_{0}, \\
               & u_{k} \in \mathcal{U},\, \forall k \in [0, \dots, M-1], \\
               & Z_{n_{0},k} \in \mathcal{Z},\, \forall k \in [0, \dots, M-1], \\
               & Z_{n_{0},M} \in \mathcal{Z}_{f}, 
\end{align}
\noindent where $M$ is the control horizon, $\bar{P}$ is the terminal cost and $\mathcal{Z}_{f}$ is the terminal set constraint. 



Finally, to evaluate our method, we present two numerical scenarios. The training data was generated by sampling the MPC's control policy over its feasible set. Specifically, we defined a grid of initial conditions and, for each point, recorded the initial state and the corresponding optimal control action \textit{only if} the MPC could find a feasible solution. To create the dataset, we generated trajectories using the state penalty $\bar{Q} = 2I_{n_{0}}$, while the control penalty $\bar{R} = 1$.

\subsection{NN Controller for the Unstable Parabolic PDE}
Starting with a more conservative approach. Setting up $q_{c} = 24$ and $\delta = 5$, we have by the relation that the number of modes $-\lambda_{n} + q_{c} < -\delta < 0$ for all $n \geq n_{0}+1$, i.e., this results in two unstable modes $n_0 = 2$. For the MPC controller, we impose the set of states $\mathcal{Z} = \left\{ z \in \mathbb{R}^{2} ~|~ \left|z_{1}\right| \leq 2,\, \left|z_{2}\right| \leq 40 \right\}$ and control actions constraints defined as $\mathcal{U} = \left\{ u \in \mathbb{R} ~|~  \left|u\right| \leq 20 \right\}$.

\begin{figure}[htpb]
    \centering
    \includegraphics[width=0.9\linewidth]{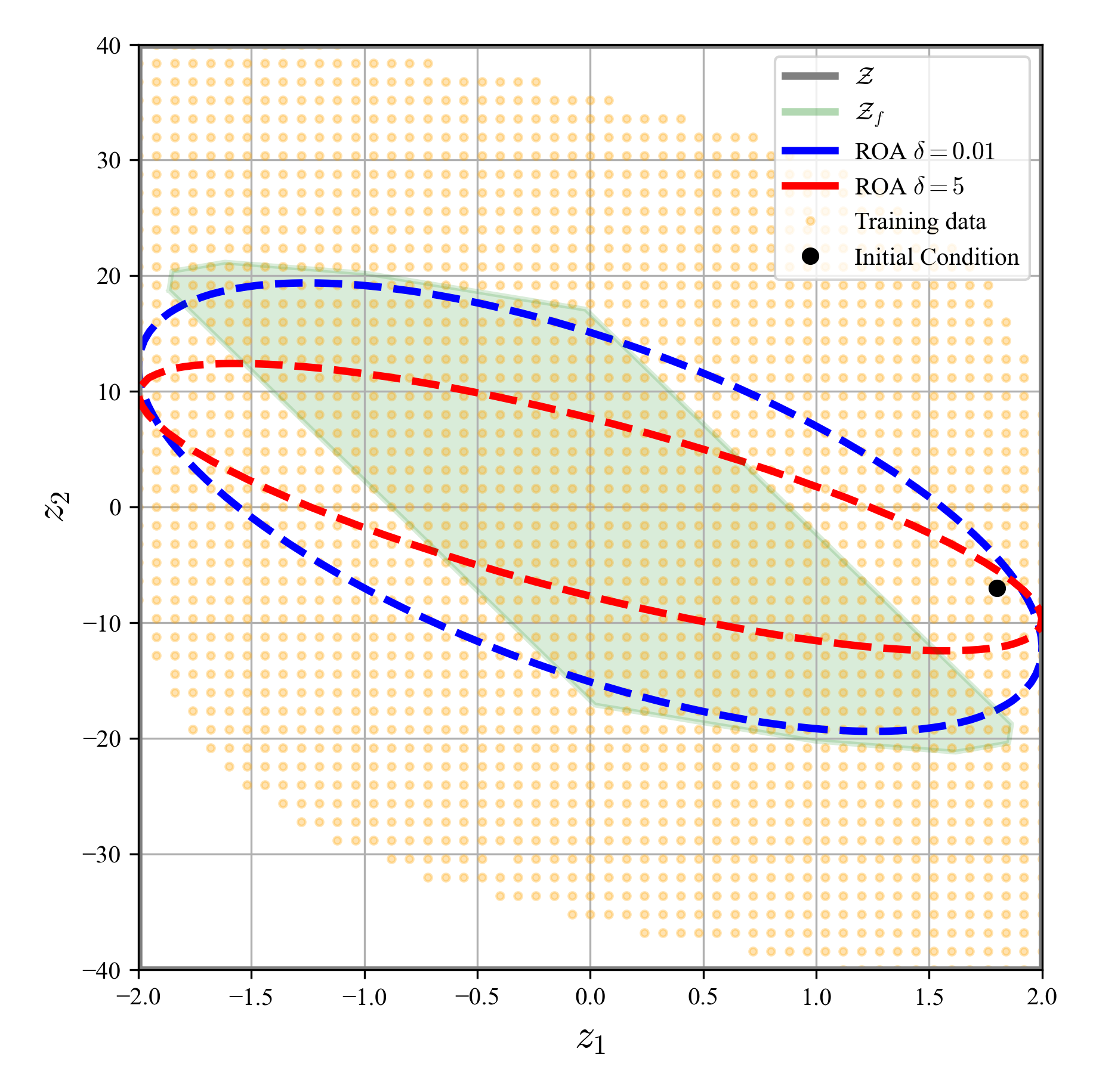}
    \caption{ROA of the NN for the first simulation example.}
    \label{fig:1}
\end{figure}

\begin{figure}[htpb]
    \centering
    \includegraphics[width=\linewidth]{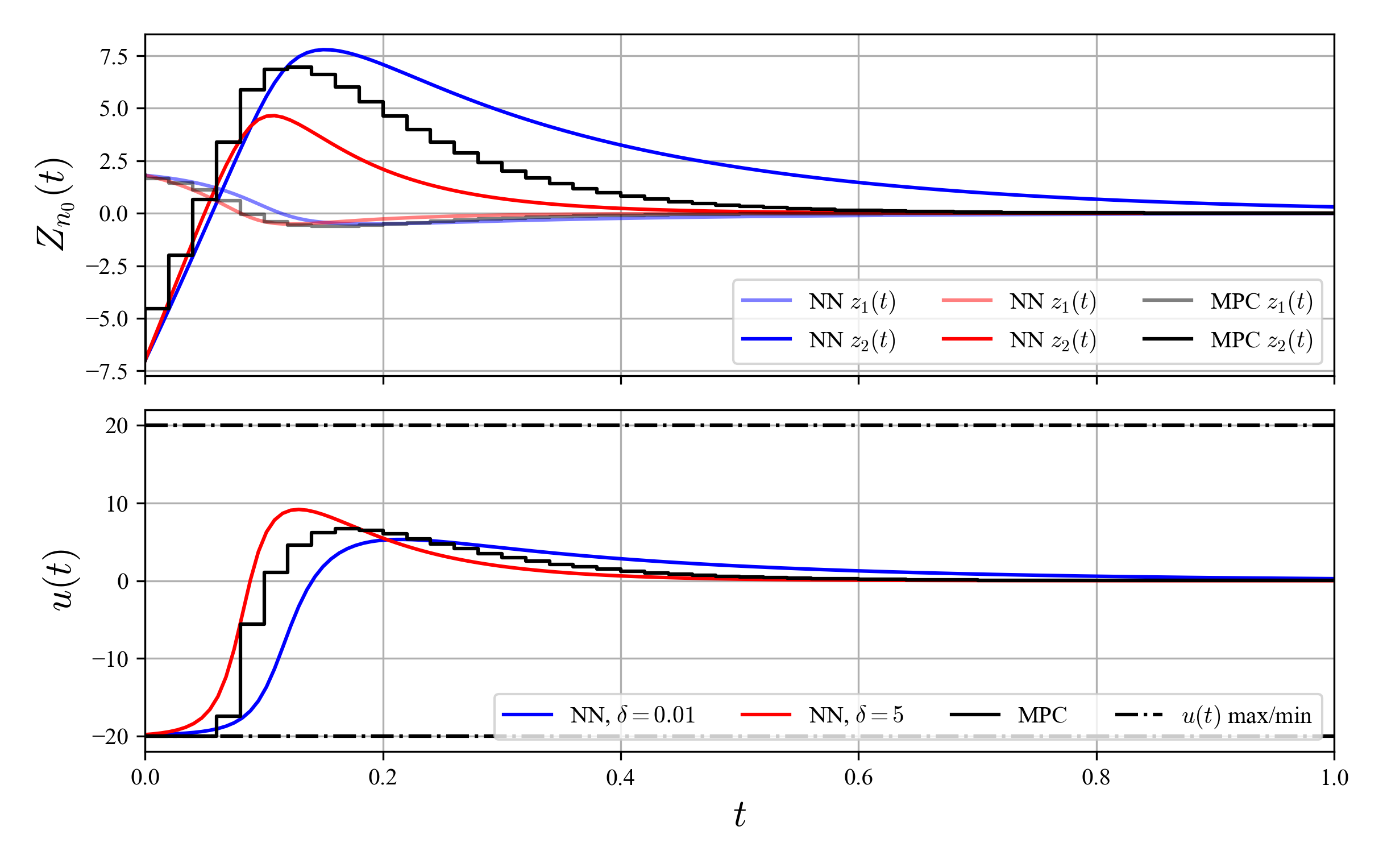}
    \caption{The dynamics of the NN and the MPC.}
    \label{fig:2}
\end{figure}


Thus, after conducting the data collection, we perform the Algorithm \ref{alg:safe_il_pde} with the tuning parameters $\eta_{1} = 1$, $\eta_{2} = 200$, $\rho = 1$, $\gamma = 1$, and $\tau = 0.1$. As a result, the graph of Figure \ref{fig:1} presents the phase state, illustrating how the ROA (blue line) lies within the dataset, the domain $\mathcal{Z}$, and the terminal set $\mathcal{Z}_f$ from the MPC. To compare the results, the MPC expert that the NN imitates is shown in black, as illustrated in Figure \ref{fig:2}. In this scenario, all LMI constraints were satisfied.





Regarding the second experiment, we aim to find a larger ROA. To achieve this, we introduce a less restrictive stability margin of $\delta = 0.01$, as can be seen in the larger ROA and the dynamic evolution (in red) in Figures \ref{fig:1} and \ref{fig:2}, respectively. However, in this scenario, we can not find a solution that satisfies the condition established at \eqref{residual:op_2}. To overcome this, we increase the number of modes by adding one stable mode, i.e., $n_0 + 1$, to relax the constraint \eqref{eq:lmi3}.

\subsection{Increasing the Number of Modes}

\begin{figure}[t]
\vspace{5pt}
    \centering
    \includegraphics[width=\linewidth]{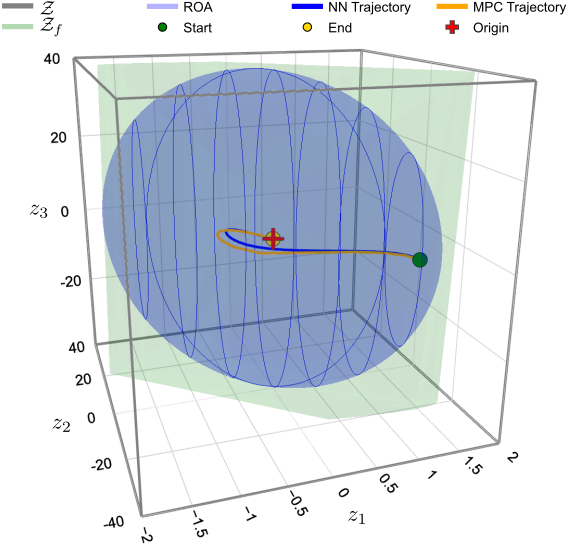}
    \caption{ROA of the NN and the dynamics of the MPC and NN for the simulation with $n_0+1$ modes.}
    \label{fig:5}
\end{figure}

In this last scenario, we increase the training process by adding one stable mode to the learning process of the NN. The idea is to keep a larger ROA while satisfying the inequality \eqref{residual:op_2}. Therefore, defining the set $\mathcal{Z}' = \mathcal{Z} \times \mathcal{S}$ and the terminal set $\mathcal{Z}'_f$, where $\mathcal{S} = \left\{ z_3 \in \mathbb{R} ~|~ \left|z_3\right| \leq 40 \right\}$, and $\mathcal{Z}' \subset \mathbb{R}^{3}$. Also, we set up the stability margin $\delta=0.1$ and the parameters $\eta_1 = 50$, $\eta_2=100$, and $\rho = 10$.

\begin{figure}[t]
\vspace{1pt}
    \centering
    \includegraphics[width=\linewidth]{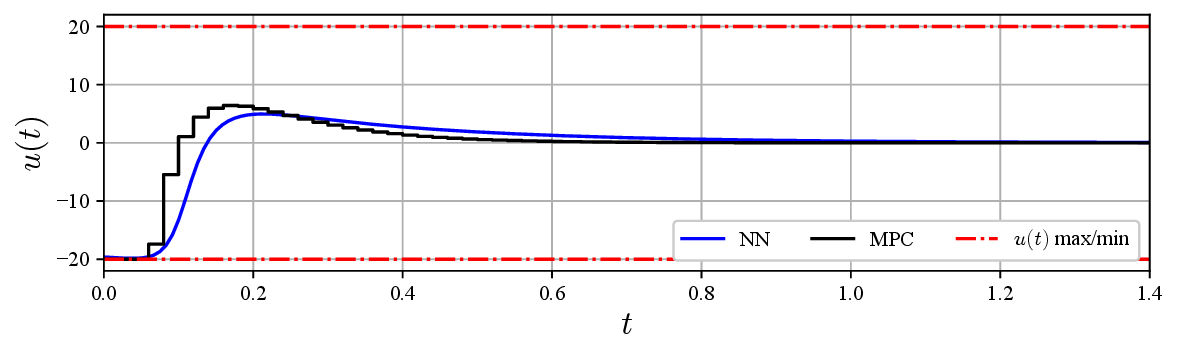}
    \caption{NN controller with $\delta=0.1$.}
    \label{fig:6}
\end{figure}

The inequalities are satisfied, and Figure \ref{fig:5} presents a larger ROA and the dynamics evolving in time. Figure \ref{fig:6} shows the control actions.

We perform simulations to evaluate the computational time requirements for each method. The NN has two layers with 10 neurons each, and for the MPC, we set up a horizon $M = 50$. The results are stated in Table \ref{tab:timing}. 
\begin{table}[h!]
\centering
\caption{Controller Execution Time.}
\label{tab:timing}
\begin{tabular}{l c c}
\toprule
Number of Modes & ($n_0=2$) & ($n_0+1$) \\
\midrule
MPC Avg. Time (ms) & 11.13 & 13.73 \\
NN Avg. Time (ms)  & 0.011 & 0.012 \\
\midrule
\textbf{Speed-up Factor} & \textbf{1024x} & \textbf{1143x} \\
\bottomrule
\end{tabular}
\end{table}

\section{Conclusion}
In this paper, we have introduced a framework for designing NN controllers for systems governed by PDE, providing formal, continuous-time stability guarantees. By integrating Lyapunov-based stability conditions directly into the learning process, we successfully synthesized a controller that is provably stable within a verifiable ROA. Crucially, our theoretical results address the challenge of model truncation, ensuring that the learned policy safely stabilizes the true infinite-dimensional system, not just its finite-order approximation. 

From the results of this work, new questions arise that open paths for new research.  A primary issue is to consolidate the separate LMI conditions for stability/robustness and the influence of the neglected high-order modes into a unified formulation, which could potentially yield less conservative guarantees. Furthermore, the framework can be extended to tackle more complex PDE systems with disturbances and uncertainties, or even observers. 



\bibliography{bibli}

\end{document}